\documentclass[12pt]{article}

\usepackage{amsmath,amssymb,amsthm}
\usepackage {graphicx}

\newtheorem{theorem}{Theorem}[section]
\newtheorem{lemma}[theorem]{Lemma}
\newtheorem{remark}[theorem]{Remark}

\newtheorem {proposition}[theorem]{Proposition}
\newtheorem {definition}[theorem]{Definition}
\newtheorem{problem}[theorem]{Problem}
\def \RR{\mathbb R}
\def \NN{\mathbb N}
\def\({\left(}  \def\){\right)}
\def\[{\left[}  \def\]{\right]}

\def \beq {\begin {equation}}
\def \eeq {\end{equation}}
\def \OL {\overline}

\def \W {\widetilde}

\begin {document}
\large{

\title {Universal minima of potentials of certain spherical designs contained in the fewest parallel hyperplanes}

\maketitle

\begin {center}
\author {Sergiy Borodachov}

\bigskip

{\it Department of Mathematics, Towson University, Towson, MD, 21252}

\abstract {We find the set of all universal minimum points of the potential of the $16$-point sharp code on $S^4$ and (more generally) of the demihypercube on $S^d$, $d\geq 5$, as well as of the $2_{41}$ polytope on $S^7$. We also extend known results on universal minima of three sharp configurations on $S^{20}$ and $S^{21}$ containing no antipodal pair to their symmetrizations about the origin. Finally, we prove certain general properties of spherical $(2m-1)$-designs contained in as few as $m$ parallel hyperplanes (all but one configuration considered here possess this property).}

\end {center}



\section {Statement of the problem and review of known results}

Let $S^d:=\{(x_1,\ldots,x_{d+1})\in \RR^{d+1}: x_1^2+\ldots+x_{d+1}^2=1\}$ be the unit sphere in $\RR^{d+1}$. Throughout the text, $\omega_N$ will denote a configuration of $N$ pairwise distinct points on $S^d$ and ${\bf x}_1,\ldots,{\bf x}_N$ will denote the points in $\omega_N$. We will also call $\omega_N$ a spherical code.
We call a function $g:[-1,1]\to(-\infty,\infty]$ an {\it admissible potential function} if $g$ is continuous on $[-1,1)$ with $g(1)=\lim\limits_{t\to 1^-}g(t)$ and differentiable in $(-1,1)$. Additional assumption(s) on the derivative(s) of $g$ will be further specified in each case. Define the $g$-potential of $\omega_N$ by
$$
p^g({\bf x},\omega_N):=\sum\limits_{i=1}^{N}g({\bf x}\cdot{\bf x}_i),\ \ \ {\bf x}\in S^d.
$$
We consider the following extremal problem over the sphere.
\begin {problem}\label {P1}
{\rm
Find the quantity
\begin {equation}\label {1}
P^g(\omega_N,S^d):=\min\limits_{{\bf x}\in S^d}p^g({\bf x},\omega_N)
\end {equation}
and points ${\bf x}^\ast\in S^d$ attaining the minimum in \eqref {1}.
}
\end {problem}

We have an important special case when the potential function $g$ is absolutely monotone on $[-1,1)$. Then $g(t)=f(2-2t)$ for some completely monotone $f$ on $(0,4]$. Recall that a function $g$ is called {\it absolutely monotoe} or {\it completely monotone} on an interval $I$ if $g^{(k)}\geq 0$ or $(-1)^kg^{(k)}\geq 0$ on $I$, respectively, for every $k\geq 0$. The function $g$ is {\it strictly absolutely} or {\it strictly completely monotone} on $I$ if the corresponding inequality is strict in the interior of $I$ for all $k\geq 0$. The kernel 
\begin {equation}\label {*}
g({\bf x}\cdot{\bf y})=f\(\left|{\bf x}-{\bf y}\right|^2\)
\end {equation}
with a strictly absolutely monotone $g$ on $[-1,1)$ includes the Riesz $s$-kernel for $s>0$ and the Gaussian kernel. After adding an appropriate positive constant, it also includes the logarithmic kernel and the Riesz $s$-kernel for $-2<s<0$.

The problem about absolute extrema on the sphere of potentials of spherical codes was earlier solved by Stolarsky \cite {Sto1975circle,Sto1975} and Nikolov and Rafailov \cite {NikRaf2011,NikRaf2013} for Riesz $s$-kernels, $s\neq 0$, and sets of vertices of a regular $N$-gon on $S^1$ and of a regular simplex, regular cross-polytope, and cube inscribed in $S^d$.
Hardin, Kendall, and Saff \cite {HarKenSaf2013} proved that absolute minima of the potential of a regular $N$-gon on $S^1$ 
with respect to a decreasing and convex function of the geodesic distance are attained at points of the dual regular $N$-gon. 

Recently, results from \cite {Sto1975circle,Sto1975,NikRaf2011,NikRaf2013} for $s>-2$, $s\neq 0$ (except for absolute maxima for the cube), were extended to kernels \eqref {*} with an absolutely monotone function $g$ and spherical designs of the highest (in a certain sense) strength. In particular, for a regular simplex on $S^d$, absolute maxima are at its vertices and absolute minima are at the antipods of its vertices, see \cite {Borsimplex}. Absolute maxima with respect to kernel \eqref {*} were found in \cite {BorMinMax} for sharp spherical codes that are antipodal or are designs of an even strength (called by some authors ``strongly sharp"). Absolute maxima appear to be independent of the potential function $g$ when $g$ is strictly absolutely monotone on $[-1,1)$ (they are at points of the code itself). Such absolute maxima are called {\it universal maxima}. 

The set of universal minima (defined in a similar way) of any spherical code that, for some $m\in \NN$, is a $(2m-1)$-design forming $m$ distinct dot products with some point ${\bf z}\in S^d$, is also known (we will call such codes $m$-stiff). It is exactly the set of all such points ${\bf z}$, see talk \cite {Bor2022talk}\footnote {Talk \cite {Bor2022talk} was given in January, 2022 at ESI and can be found in the ESI's YouTube accound.} by the author for the proof or Lemma~3.5 in \cite {BorMinMax}\footnote{Paper \cite {BorMinMax} is on ArXiv since March 2022.}, from which this result follows.  We will call the set of such points ${\bf z}$ the dual configuration. We restate this result here as Theorem~\ref {2m-1}. Immediate consequences of this result are universal minima of a regular $2m$-gon on $S^1$, of a regular cross-polytope and cube on $S^d$, and of the $24$-cell on $S^3$, since finding the dual configuration is elementary for these codes. Stiff spherical codes attain the Fazekas-Levenshtein bound for covering \cite [Theorem 2]{FazLev1995}.

Universal minima of any strongly sharp code are at the antipods of points of the code. This  also follows from \cite [Lemma 3.5]{BorMinMax}. Immediate consequences of this fact (other than the regualr simplex) are universal minima of a regular $(2m+1)$-gon on $S^1$, the Schl\"affi configuration on $S^5$, and the McLaughlin configuration on~$S^{21}$. 
Strongly sharp spherical codes also attain the Fazekas-Levenshtein bound for covering \cite [Theorem 2]{FazLev1995}.

Boyvalenkov, Dragnev, Hardin, Saff, and Stoyanova \cite [Theorems 3.4 and 3.7]{BoyDraHarSafSto600cell} (see also \cite [Theorem 1.4]{BoyDraHarSafStosharpantipodal})\footnote {Papers \cite {BoyDraHarSafSto600cell} and \cite {BoyDraHarSafStosharpantipodal} are on ArXiv since July and October 2022, respectively.} proved universal upper and lower bounds for the potential of a general spherical design. These bounds become sharp in the cases mentioned above: in the case of a minimum for stiff and strongly sharp configurations and in the case of a maximum for sharp antipodal and strongly sharp ones. The lower bound is an analogue of the Fazekas-Levenshtein bound for covering \cite[Theorem 2] {FazLev1995}.

Paper \cite {BoyDraHarSafSto600cell} also showed that the universal maxima of the $600$-cell on $S^3$ are vertices of the $600$-cell itself. The work \cite {BoyDraHarSafStosharpantipodal} proved that a number of known sharp codes are also stiff. Then the result from \cite {Bor2022talk} implies that their universal minima are at points of the dual configuration. Paper \cite {BoyDraHarSafStosharpantipodal} further studies the dual configuration for each code and antipods of the two strongly sharp codes on $S^5$ and $S^{21}$ mentioned above. The author in paper \cite {Borstiff}\footnote {On ArXiv since December 2022.} found explicitly the sets of all universal minima for five more stiff configurations (which are not sharp) on spheres of different dimensions as well as for the $56$-point kissing configuration on $S^6$, which is a known sharp code (paper \cite {BoyDraHarSafStosharpantipodal} gives one universal minimum of this code).

Certain remarkable spherical configurations are not stiff or strongly sharp. Universal minima of the regular icosahedron and regular dodecahedron on $S^2$ were characterized in \cite {Borsymmetric}\footnote {On ArXiv since October 9, 2022. The proof for icosahedron was briefly discussed in talk \cite {Bor2022talk} in January 2022.} as well as universal minima of the $E_8$ lattice on $S^7$. Furthermore, one universal minimum and the corresponding absolute minimum value of the potential were found in \cite {BoyDraHarSafStosharpantipodal}\footnote {On Arxiv since October 31, 2022.} for the Leech lattice on $S^{23}$. Papers \cite {Borsymmetric}\footnote {The general ``skip one add two" theorem from \cite {Borsymmetric} was stated in talk \cite {BorBatumi} in August 2022, see Theorem \ref {generalth}.} and \cite {BoyDraHarSafStosharpantipodal} establish general theorems (diferent to a certain extent) for the so-called ``skip one add two" case and use them to establish the above mentioned results. 

Critical points of the total potential of finite configurations of charges were also analysed (see \cite {Bil2015,GioKhi2020} and references therein). This work is related to the known Maxwell's conjecture.
A more detailed review of known results on extrema of potentials of spherical codes can be found, for example, in \cite {Borstiff}.

In this paper, we prove certain properties of general stiff configurations and characterize universal minima of the $16$-point sharp code\footnote {In \cite {BoyDraHarSafStosharpantipodal}, its universal minima were presented without a characterization.} on $S^4$, of a demihypercube on $S^d$, $d\geq 5$, and of the $2_{41}$ polytope on $S^7$ (it is dual to the $E_8$ lattice). For some sharp spherical codes that have no antipodal pairs, we show that universal minima found in \cite {BoyDraHarSafStosharpantipodal} are also universal minima for the symmetrization of each of them about the origin.

One of important applications of Problem \ref {P1} is the polarization problem on the sphere. Papers \cite {Sto1975circle,NikRaf2011,HarKenSaf2013,Borsimplex,BorMinMax,BoyDraHarSafSto600cell} that we mentioned when reviewing results on extrema of potentials solve certain its cases. A more comprehensive review of known work on polarization can be found, for example, in book \cite [Chapter~14] {BorHarSafbook} with most recent results reviewed in, e.g., \cite {BorMinMax}.

The paper is structured as follows. Section \ref {prel} contains the necessary preliminaries. In Section \ref {cube}, we characterize universal minima of the $d$-demicube for $d\geq 5$ (including the $16$-point sharp code on $S^4$). In Section \ref {2_41}, we find all the universal minima of the $2_{41}$ polytope on~$S^7$. Section~\ref {sy} extends known results on universal minima for certain three non-antipodal sharp configurations to their symmetrizations. In Section \ref {scstiff}, we establish certain properties of general stiff configurations and of their duals.

\section {Preliminaries}\label {prel}

In this section, we state definitions and known facts used further in the paper. Define 
$$
w_{d}(t):=\gamma_{d}(1-t^2)^{d/2-1}, 
$$
where the constant $\gamma_{d}$ is such that $w_{d}$ is a probability density on $[-1,1]$. The {\it Gegenbauer orthogonal polynomials} corresponding to the sphere $S^{d}$ in $\RR^{d+1}$ are terms of the sequence $\{P_n^{(d)}\}_{n=0}^{\infty}$ of univariate polynomials such that ${\rm deg}\ \! P_n^{(d)}=n$, $n\geq 0$, and 
$$
\int_{-1}^{1}P_i^{(d)}(t)P_j^{(d)}(t)w_{d}(t)\ \! dt =0,\ \ \ i\neq j,
$$
normalized so that $P_n^{(d)}(1)=1$, $n\geq 0$ (see \cite [Chapter 4]{Sze1975} or, e.g., \cite [Chapter~5]{BorHarSafbook}). 

For a configuration $\omega_N=\{{\bf x}_1,\ldots,{\bf x}_N\}\subset S^d$, let $\mathcal I(\omega_N)$ be the set of all $n\in \NN$ such that
\begin {equation}\label {Pge}
\sum\limits_{i=1}^{N}\sum\limits_{j=1}^{N}P_n^{(d)}({\bf x}_i\cdot {\bf x}_j)=0.
\end {equation}
We call $\mathcal I(\omega_N)$ {\it the index set} of $\omega_N$.
Let $\sigma_{d}$ be the $d$-dimensional area measure on the sphere $S^d$ normalized to be a probability measure.
A configuration $\omega_N$ is called a {\it spherical $n$-design} if, for every polynomial $p$ on $\RR^{d+1}$ of degree at most~$n$,
\begin {equation}\label {design}
\frac {1}{N}\sum\limits_{i=1}^{N}p({\bf x}_i)=\int_{S^{d}}p({\bf x})\ \! d\sigma_{d}({\bf x}),
\end {equation}
see the paper by Delsarte, Goethals, and Seidel \cite {DelGoeSei1977}. The maximal number $n$ in this definition is called the {\it strength} of the spherical design $\omega_N$.

We recall the following equivalent definitions of a spherical design. Let $\mathbb P_n$ denote the space of all univariate polynomials of degree at most $n$. 
\begin {theorem}\label {cd}(see \cite {DelGoeSei1977,FazLev1995} or, e.g., \cite[Lemma 5.2.2 and Theorem 5.4.2]{BorHarSafbook})
Let $d,n\geq 1$ and $\omega_N=\{{\bf x}_1,\ldots,{\bf x}_N\}$ be a point configuration on $S^d$. The following are equivalent:
\begin {itemize}
\item[(i)]
$\omega_N$ is a spherical $n$-design;
\item [(ii)]
$
\{1,\ldots,n\}\subset \mathcal I(\omega_N);
$
\item[(iv)]
for every polynomial $q\in \mathbb P_n$, we have $p^q({\bf y},\omega_N)=\sum_{i=1}^{N}q({\bf y}\cdot{\bf x}_i)=C$, ${\bf y}\in S^d$, where $C$ is a constant.
\end {itemize}
\end {theorem}
If item (iii) holds in the above theorem, then $C=a_0(q)N$, where
\begin {equation}\label {a0}
a_0(q):=\int_{-1}^{1}q(t)w_d(t)\ \! dt
\end {equation}
is the $0$-th Gegenbauer coefficient of polynomial $q$. For a given $m\in \NN$ and a given configuration $\omega_N\subset S^{d}$,
denote by $\mathcal D_m(\omega_N)$ the set of all points ${\bf z}\in S^{d}$ for which the set of dot products 
$$
D({\bf z},\omega_N):=\{{\bf z}\cdot {\bf x}_i : i=1,\ldots,N\}
$$ 
has at most $m$ distinct elements. 
\begin {definition}\label {D1}
{\rm
We call a point configuration $\omega_N\subset S^{d}$ {\it $m$-stiff}, $d,m\geq 1$, if $\omega_N$ is a spherical $(2m-1)$-design and the set $\mathcal D_m(\omega_N)$ is non-empty. The set $\mathcal D_m(\omega_N)$ of a given $m$-stiff configuration $\omega_N$ is called {\it the dual configuration} for~$\omega_N$.
}
\end {definition}
Following \cite {CohKum2007}, we call a configuration $\omega_N\subset S^d$ {\it sharp} if, for some $m\in \mathbb N$, it is a $(2m-1)$-design and there are exactly $m$ distinct values of the dot product between distinct points in $\omega_N$. If, in addition, $\omega_N$ is a $2m$-design, we call it {\it strongly sharp}.

The next statement is a part of the classification of quadratures in \cite {DelGoeSei1977}  corresponding to spherical designs of the highest strength; i.e., stiff, strongly sharp, or sharp antipodal codes. For the reader's convenience, we mention that its proof can also be found, for example, in \cite [Proposition~7.2]{Borstiff}.
Let $\{\varphi_1,\ldots,\varphi_m\}$ be the fundamental polynomials for the set $-1<\kappa_1^{m}<\ldots<\kappa_{m}^{m}<1$ of zeros of the Gegenbauer polynomial $P_m^{(d)}$; that is, $\varphi_i\in \mathbb P_{m-1}$, $\varphi_i(\kappa_i^m)=1$, and $\varphi_i(\kappa_j^m)=0$, $j\neq i$, $i=1,\ldots,m$.
\begin {proposition}\label {mtight}
If $\omega_N$ is an $m$-stiff configuration on $S^d$, then for every ${\bf z}\in \mathcal D_m(\omega_N)$, the set $D({\bf z},\omega_N)$ contains exactly $m$ distinct elements located in $(-1,1)$, which are $\kappa_1^{m},\ldots,\kappa_{m}^{m}$. Furthermore, the number of indices $i$ such that ${\bf z}\cdot {\bf x}_i=\kappa_j^{m}$ does not depend on ${\bf z}$ and equals $a_0(\varphi_j) N$, $j=1,\ldots,m$.

In particular, if $m=2$, then for every ${\bf z}\in \mathcal D_2(\omega_N)$, we have $D({\bf z},\omega_N)=\left\{-\frac {1}{\sqrt{d+1}},\frac {1}{\sqrt{d+1}}\right\}$.
\end {proposition}

\begin {remark}
{\rm 
In view of Proposition \ref {mtight}, an $m$-stiff configuration may exist on $S^d$ for given $m,d\geq 1$, only if the all the numbers $a_0(\varphi_i)$, $i=1,\ldots,m$, are positive rationals. These numbers are the weights of the Gauss-Gegenbauer quadrature for integral \eqref {a0} (the nodes are $\kappa_1^m,\ldots,\kappa_m^m$).
}
\end {remark}

We next restate the result proved in talk \cite {Bor2022talk}. It follows from \cite [Lemma 3.5]{BorMinMax}, see the proof of \cite [Theorem~4.3]{Borstiff} for details.
\begin {theorem}\label {2m-1}
Let $m\geq 1$, $d\geq 1$, and $g$ be an admissible potential function with a convex derivative $g^{(2m-2)}$ on $(-1,1)$. If $\omega_N=\{{\bf x}_1,\ldots,{\bf x}_N\}$ is an $m$-stiff configuration on the sphere $S^{d}$, then the potential
$$
p^g({\bf x},\omega_N)=\sum\limits_{i=1}^{N}g({\bf x}\cdot{\bf x}_i),\ \ \ {\bf x}\in S^{d}, 
$$
attains its absolute minimum over $S^{d}$ at every point of the set $\mathcal D_m(\omega_N)$. 

If, in addition, $g^{(2m-2)}$ is strictly convex on $(-1,1)$, then $\mathcal D_m(\omega_N)$ contains all points of absolute minimum of the potential $p^g(\cdot,\omega_N)$ on $S^{d}$.
\end {theorem}

We will also need the ``skip one add two" result from \cite [Theorem 3.1]{Borsymmetric}. 

\begin {theorem}\label {generalth}
Let $d,m\in \NN$, $m\geq 2$, and $\omega_N=\{{\bf x}_1,\ldots,{\bf x}_N\}$ be a point configuration on $S^d$ whose index set $\mathcal I(\omega_N)$ contains numbers $1,2,\ldots,2m-3,2m-1,2m$. Assume that numbers $-1< t_1<t_2<\ldots<t_m< 1$ are such that 
\begin {equation}\label {sumsq}
\sum_{i=1}^{m}t_i<t_m/2 \ \ \text { and }\ \ \sum\limits_{i=1}^{m}t_i^2-2\(\sum\limits_{i=1}^{m}t_i\)^2<\frac {m(2m-1)}{4m+d-3},
\end {equation}
and that the set $\mathcal D$ of points ${\bf x}^\ast\in S^d$ with $D({\bf x}^\ast,\omega_N)\subset \{t_1,\ldots,t_m\}$ is non-empty.
Let $g$ be an admissible potential function with non-negative derivatives $g^{(2m-2)}$, $g^{(2m-1)}$, and $g^{(2m)}$ on $(-1,1)$. Then, for every point ${\bf x}^\ast\in \mathcal D$,
\begin {equation}\label {star}
\min\limits_{{\bf x}\in S^d}\sum\limits_{i=1}^{N}g({\bf x}\cdot{\bf x}_i)=\sum\limits_{i=1}^{N}g({\bf x}^\ast\cdot {\bf x}_i).
\end {equation}
If, in addition, $g^{(2m)}>0$ on $(-1,1)$, then the absolute minimum in \eqref {star} is achieved only at points of the set $\mathcal D$. 
\end {theorem}

We remark that proofs of Theorems \ref {2m-1} and \ref {generalth} utilize the Delsarte-Yudin method (also known as the Delsarte or linear programming or polynomial method), see the work by Delsarte, Goethals, and Seidel \cite {DelGoeSei1977} or by Yudin \cite {Yud1992}. A detailed description of this approach and references to works using it can also be found, in particular, in \cite {Lev1979,Lev1992,Lev1998,CohKum2007,BoyDraHarSafSto600cell,BoyDraHarSafStosharpantipodal} and in \cite [Chapter 5]{BorHarSafbook}.

\section {The $16$-point sharp code on $S^4$ and the demihypercube}\label {cube}

Denote by $\omega_{2d}^\ast:=\{\pm {\bf e}_1,\ldots,\pm {\bf e}_d\}$, $d\geq 2$, where ${\bf e}_1,\ldots,{\bf e}_d$ are vectors of the standard basis in $\RR^d$, the set of vertices of the regular cross-polytope inscribed in $S^{d-1}$ and let $U_d$ be the set $\left\{\(\pm \frac {1}{\sqrt {d}},\ldots,\pm \frac {1}{\sqrt {d}}\)\right\}\subset \RR^d$ of vertices of the cube inscribed in $S^{d-1}$. It is not difficult to see that $\mathcal D_2(\omega_{2d}^\ast)=U_d$ and $\mathcal D_2(U_d)=\omega_{2d}^\ast$.

Let $\OL\omega^d$, $d\geq 2$, be the set of $N=2^{d-1}$ points $\(\pm \frac {1}{\sqrt{d}},\ldots,\pm \frac {1}{\sqrt{d}}\)\in U_d$ with an even number of minus signs. This configuration forms the set of vertices of a {\it $d$-demicube} (also called the {\it demihypercube}). The set $\W\omega^d$ of vectors from $U_d$ with an odd number of minus signs is a reflection of $\OL\omega^d$ with respect to any of the coordinate hyperplanes; i.e., $\W\omega^d$ is an isometric copy of $\OL\omega^d$. Therefore, it is sufficient to consider just $\OL\omega^d$. We have $\OL\omega^d\cup \W\omega^d=U_d$ and the two sets are disjoint. For $d$ odd, we have $\W\omega^d=-\OL\omega^d$ with both sets not containing antipodal pairs. For $d$ even, each of the sets $\OL\omega^d$ and $\W\omega^d$ is itself antipodal.

Observe that for $d=2$, both configurations consist of one antipodal pair; i.e., they are $1$-stiff. For $d=3$, each one is a regular simplex inscribed in $S^2$ (strongly sharp and, hence, not stiff).
For $d=4$, the set $\OL\omega^d$ consists of eight points $\(\pm \frac 12,\pm \frac 12,\pm \frac 12,\pm \frac 12 \)$ with an even number of minus signs. Each set $\OL\omega^4$ and $\W\omega^4$ is an isometric copy of a regular cross-polytope in $\RR^4$; i.e., it is $2$-stiff. 
For $d=5$, configuration $\OL\omega^d$ consists of $16$ points on $S^4$ of the form $\(\pm \frac {1}{\sqrt{5}},\pm \frac {1}{\sqrt{5}},\pm \frac {1}{\sqrt{5}},\pm \frac {1}{\sqrt{5}},\pm \frac {1}{\sqrt{5}}\)$ with an even number of minus signs. This is the well-known sharp $(5,16,1/5)$-code. It was described by Gossett \cite {Gos1900}. The set $\W\omega^5$ is the antipode of this code. The $2$-stiffness property of $\OL\omega^5$ was observed in \cite {BoyDraHarSafStosharpantipodal}.
We show that the $d$-demicube is $2$-stiff for any $d\geq 6$. We start with the following auxiliary statement.
\begin {lemma}\label {3design}
Let $\omega_N$ be a non-empty subset of $U_d=\left\{-\frac {1}{\sqrt{d}},\frac {1}{\sqrt{d}}\right\}^d$, $d\geq 3$. Then $\omega_N$ is a $3$-design if, and only if, $N$ is even and for every set $I$ of one, two, or three pairwise distinct indices, exactly half of vectors in $\omega_N$ have an even number of negative coordinates with indices in $I$ and exactly half have an odd number of negative 
coordinates with indices in $I$.

\end {lemma}
\begin {proof}
Let $\omega_N\subset U_d$ be arbitrary.
Using the notation ${\bf y}=(y_1,\ldots,y_d)$ for a point ${\bf y}\in \omega_N$, define
$$
S_i:=\sum\limits_{{\bf y}\in \omega_N}y_i,\ \ \ \ S_{i,j}:=\sum\limits_{{\bf y}\in \omega_N}y_iy_j, \ \ \ \text{and}\ \ \ S_{i,j,k}:=\sum\limits_{{\bf y}\in \omega_N}y_iy_jy_k.
$$
Observe that $S_{i,j}$ and $S_{i,j,k}$ do not depend on permutations of indices and that $S_{i,i}=\frac {N}{d}$. When some two indices coincide, say $i=j$, we have 
\begin {equation}\label {i=j}
S_{i,j,k}=\sum\limits_{{\bf y}\in \omega_N}\frac {y_k}{d}=\frac {1}{d}S_k.
\end {equation}
Formula \eqref {i=j} holds even if $i=j=k$.
Let ${\bf x}=(x_1,\ldots,x_d)\in S^{d-1}$ be any vector. Then
$$
\sum\limits_{{\bf y}\in \omega_N}{\bf x}\cdot {\bf y}=\sum\limits_{{\bf y}\in \omega_N}\sum\limits_{i=1}^{d}x_iy_i=\sum\limits_{i=1}^{d}\sum\limits_{{\bf y}\in \omega_N}x_iy_i=\sum\limits_{i=1}^{d}x_i\sum\limits_{{\bf y}\in \omega_N}y_i=\sum\limits_{i=1}^{d}S_ix_i,
$$
\begin {equation*}
\begin {split}
\sum\limits_{{\bf y}\in \omega_N}&\({\bf x}\cdot {\bf y}\)^2=\sum\limits_{{\bf y}\in \omega_N}\(\sum\limits_{j=1}^{d}x_jy_j\)^2=\sum\limits_{{\bf y}\in \omega_N}\sum\limits_{i=1}^{d}\sum\limits_{j=1}
^{d}x_ix_jy_iy_j=\sum\limits_{i=1}^{d}\sum\limits_{j=1}
^{d}\sum\limits_{{\bf y}\in \omega_N}x_ix_jy_iy_j\\
&=\sum\limits_{i=1}^{d}x_i^2\sum\limits_{{\bf y}\in \omega_N}y_i^2+\sum\limits_{i,j=1\atop i\neq j}^{d}x_ix_j\sum\limits_{{\bf y}\in \omega_N}
y_iy_j=\frac {N}{d}+\sum\limits_{i,j=1\atop i\neq j}^{d}S_{i,j}x_ix_j,
\end {split}
\end {equation*}
and
\begin {equation*}
\begin {split}
\sum\limits_{{\bf y}\in \omega_N}\({\bf x}\cdot {\bf y}\)^3&=\sum\limits_{{\bf y}\in \omega_N}\(\sum\limits_{j=1}^{d}x_jy_j\)^3=\sum\limits_{{\bf y}\in \omega_N}\sum\limits_{i,j,k=1}^{d}x_ix_jx_ky_iy_jy_k
\\
&=\sum\limits_{i,j,k=1}^{d}x_ix_jx_k\sum\limits_{{\bf y}\in \omega_N}y_iy_jy_k=\sum\limits_{i,j,k=1}^{d}S_{i,j,k}x_ix_jx_k.
\end {split}
\end {equation*}
The configuration $\omega_N$ will be a $3$-design if, and only if, the three sums above are constant, see Theorem \ref {cd}.
If 
\begin {equation}\label {Sj}
S_i=0 \ \text{for all }i,\ \   S_{i,j}=0, \ \text{for }i\neq j, \ \ \text{and}\ \  S_{i,j,k}=0, \ \text{for}\ i,j,k \ \text {distinct}
\end {equation}
then three sums above will be constant (one should also use \eqref {i=j}). Conversely, if all three sums above are constant, then the first one has the same value for every vector $\pm {\bf e}_i$, which is $\pm S_i$, $i=1,\ldots,d$. Then $S_i=0$ for all $i$. For every vector ${\bf x}\in S^{d-1}$ with the $\ell$-th coordinate with $1/\sqrt{2}$, the $n$-th coordinate with $\pm \frac {1}{\sqrt{2}}$, and the remaining coordinates being zero, $\ell\neq n$, the value of the second sum is $N/d\pm S_{\ell,n}={\rm const}$. This forces $S_{\ell,n}=0$, $\ell\neq n$. For vector $\pm {\bf x}$, where the $\ell$-th, $n$-th, and $m$-th coordinate of ${\bf x}$ equals $1/\sqrt{3}$, $\ell,n,m$ are pairwise distinct, and the remaining coordinates are zero, the third sum equals (use \eqref {i=j} and the fact that $S_i=0$ for all $i$) 
$$
\pm \frac {1}{3\sqrt{3}}\sum_{i,j,k\in \{\ell,n,m\}}S_{i,j,k}=\pm \frac {6}{3\sqrt{3}}S_{\ell,n,m}={\rm const}. 
$$
Then $S_{\ell,n,m}=0$. Thus, $\omega_N$ is a $3$-design if, and only if, relations \eqref {Sj} hold. 

In each sum $S_i$, $S_{i,j}$, and $S_{i,j,k}$ in \eqref {Sj}, all terms have the same absolute values. Then the value of each sum in \eqref {Sj} equals that common absolute value times the difference between the number of positive and negative terms. Therefore, relations \eqref {Sj} hold if, and only if, each sum in \eqref {Sj} has equal number of positive and negative terms.
This, in turn, will hold if, and only if, for any set of indices $I=\{i\}$ or $\{i,j\}$, where $i\neq j$, or $\{i,j,k\}$, where $i,j,k$ are pairwise distinct, the number of vectors in $\omega_N$ with an even number of negative components with indices in $I$ equals the number of vectors in $\omega_N$ with an odd number of negative components with indices in~$I$. This also forces $N$ to be even.
 \end {proof}

\begin {lemma}\label {2stiff}
The $d$-demicube $\OL\omega^d$, $d\geq 4$, is $2$-stiff.
\end {lemma}
\begin {proof}
Since $\OL\omega^d$ is a subset of the set of vertices of a cube, it is contained in two parallel hyperplanes. Thus, it remains to show that $\OL\omega^d$ is a $3$-design. 
Let $I$ be any set of $k$ pairwise distinct indices, where $k=1,2,3$. A combination of signs of coordinates corresponding to $I$ with an even number of negative ones can be chosen in $2^{k-1}$ different ways. For each of these combinations, the remaining $d-k$ positions can have $2^{d-k}$ different combinations of signs with $2^{d-k-1}$ of them having an even number of minus signs. Then the total number of vectors in $\OL\omega^d$ with an even number of negative coordinates corresponding to $I$ will be $2^{k-1}\cdot 2^{d-k-1}=2^{d-2}$. By a similar argument, the number of vectors in $\OL\omega^d$ with an odd number of negative coordinates corresponding to $I$ will also be $2^{d-2}$. Lemma \ref {3design} now implies that $\OL\omega^d$ is a $3$-design and, hence, is $2$-stiff. 
\end {proof}

We next find the dual configuration of the $d$-demicube. For $d=2$, the $d$-demicube is a pair of antipodal vectors, which is $1$-stiff. Its dual is the perpendicular pair of antipodal vectors. For $d=3$, the $d$-demicube is a regular simplex, which is not stiff, since it is strongly sharp. For $d=4$, the configuration $\OL\omega^d$ is a regular cross-polytope, and its dual is the corresponding cube inscribed in $S^3$. For $d\geq 5$, we have the following result.

\begin {lemma}\label {dual}
For every $d\geq 5$, we have $\mathcal D_2(\OL\omega^d)=\omega^\ast_{2d}$. 
\end {lemma}
Since $\mathcal D_2(\omega_{2d}^\ast)=U_d$, $d\geq 2$, Lemma \ref {dual} shows that the inclusion $\omega_N\subset \mathcal D_m(\mathcal D_m(\omega_N))$ can be strict for an $m$-stiff configuration $\omega_N$ with $m\geq 2$ even if $\omega_N$ is antipodal. Furthermore, for $d\geq 5$ odd, it provides another example of a non-antipodal $m$-stiff configuration with $m\geq 2$ (non-antipodal $1$-stiff configurations are easy to construct).

\begin {proof}[Proof of Lemma \ref {dual}]
Every vector $\pm {\bf e}_i\in \omega^\ast_{2d}$ forms only dot products $\frac {1}{\sqrt{d}}$ and $-\frac {1}{\sqrt{d}}$ with points from $\OL\omega^d$; i.e., it belongs to $\mathcal D_2(\OL\omega^d)$. 
Choose any ${\bf x}=(x_1,\ldots,x_d)\in \mathcal D_2(\OL\omega^d)$. Assume to the contrary that ${\bf x}$ has at least three non-zero coordinates. Let $k$ be the number of strictly negative components in ${\bf x}$. If $k$ is even, we choose a vector ${\bf z}=(z_1,\ldots,z_d)\in \OL\omega^d$ with $-\frac {1}{\sqrt{d}}$ on all positions corresponding to strictly negative components in ${\bf x}$ and $\frac {1}{\sqrt{d}}$ on all other positions. If $k$ is odd, we choose ${\bf z}\in \OL\omega^d$ with $-\frac {1}{\sqrt{d}}$ on all but one positions corresponding to strictly negative components of ${\bf x}$ and $\frac {1}{\sqrt{d}}$ on all other positions. Then dot product ${\bf x}\cdot {\bf z}=x_1z_1+\ldots+x_dz_d$ has at most one strictly negative term and at least two other strictly positive terms, which we denote by $x_iz_i$ and $x_jz_j$. Since $d\geq 5$, we can choose two disjoint pairs of positions in ${\bf z}$ one containing $z_i$ and the other one containing $z_j$ with both pairs avoiding the position corresponding to the possible negative term in ${\bf x}\cdot {\bf z}$. Changing the sign of the coordinates of ${\bf z}$ in the first pair of positions, we keep ${\bf z}$ in $\OL\omega^d$ and strictly decrease the dot product ${\bf x}\cdot{\bf z}$. Changing the sign of the coordinates of the new vector ${\bf z}$ in the second pair of positions, we keep the resulting vector in $\OL\omega^d$ and further decrease dot product ${\bf x}\cdot {\bf z}$. This shows that ${\bf x}$ forms at least three distinct dor products with points of $\OL\omega^d$ contradicting its choice. 

Therefore, ${\bf x}$ has at most two non-zero components. Assume to the contrary that ${}\bf x$ has exactly two non-zero components, say $x_\ell$ and $x_n$. Then ${\bf x}$ forms dot products $\frac {\pm x_\ell\pm x_n}{\sqrt{d}}$ with vectors from $\OL\omega^d$ and at least three of them are distinct. Thus, ${\bf x}$ has one non-zero component. Since ${\bf x}$ is on $S^{d-1}$, this component must be $\pm 1$; that is, ${\bf x}\in \omega^\ast_{2d}$. Thus, $\mathcal D_2(\OL\omega^d)=\omega^\ast_{2d}$.
\end {proof}

We are ready to characterize universal minima of the $d$-demicube for $d\geq 5$. 
\begin {theorem}\label {demicube}
Let $d\geq 5$ and $g$ be an admissible potential function with a convex derivative $g''$ on $(-1,1)$. Then the potential
$
p^g(\cdot,\OL\omega^d)
$
of $d$-demicube $\OL\omega^d$ attains its absolute minimum over $S^{d-1}$ at every point of cross-polytope $\omega_{2d}^\ast$. 

If, in addition, $g''$ is strictly convex on $(-1,1)$, then $\omega_{2d}^\ast$ contains all points of absolute minimum of the potential $p^g(\cdot,\OL\omega^d)$ on $S^{d-1}$.
\end {theorem}
In the case $d=5$, the first paragraph of Theorem \ref {demicube} follows from the results of \cite {BoyDraHarSafStosharpantipodal}.

\begin {proof}
Since $\OL\omega^d$ is $2$-stiff, by Theorem \ref {2m-1}, the potential $p^g(\cdot,\OL\omega^d)$ attains its absolute minimum over $S^{d-1}$, $d\geq 5$, at points of the set $\mathcal D_2(\OL\omega^d)$, which, by Lemma \ref {dual}, equals $\omega^\ast_{2d}$. If $g''$ is strictly convex on $(-1,1)$, then, by Theorem~\ref {2m-1}, the set $\mathcal D_2(\OL\omega^d)=\omega^\ast_{2d}$ contains all absolute minima of $p^g(\cdot,\OL\omega^d)$ over $S^{d-1}$.
\end {proof}

\section {The $2_{41}$ polytope on $S^7$}\label {2_41}

Recall that {\it the $E_8$ lattice} is the set (lattice in $\RR^8$) of vectors in $\mathbb Z^8 \cup (\mathbb Z+{1}/{2})^8$ whose coordinates sum to an even integer.
Let $\OL \omega_{240}$ be the set of minimal length non-zero vectors of the $E_8$ lattice normalized to lie on $S^7$. The configuration $\OL \omega_{240}$ consists of $4\(8 \atop 2\)=112$ vectors with $6$ zero coordinates and two coordinates with $\pm 1/\sqrt {2}$ and $2^7=128$ vectors with all eight coordinates $\pm \frac {1}{2\sqrt {2}}$ and even number of ``$-$" signs (this part is the $8$-demicube). For brevity, we will also call $\OL \omega_{240}$ the $E_8$ lattice.

The $2_{41}$ polytope on $S^7$ (the name is due to Coxeter), denoted here by $\OL \omega_{2160}$, is the set of $N=2160$ vectors on $S^7$ that includes $16\(8 \atop 4\)=1120$ vectors with $4$ zero coordinates and $4$ coordinates with $\pm 1/2$ (let us call them type I vectors), $16$ vectors with $7$ zero coordinates and one coordinate with $\pm 1$ (let us call them type II vectors), and $8\(\(8 \atop 1\)+\(8\atop 3\)+\(8 \atop 5\)+\(8\atop 7\)\)=1024$ vectors with $7$ coordinates with $\pm 1/4$, one coordinate with $\pm 3/4$, and an odd number of negative coordinates (call them type III vectors). One can verify directly that equality \eqref {Pge} holds for $d=7$ and $n\in \{1,\ldots,7,9,10\}$. Indeed, since $\OL\omega_{2160}$ is antipodal, \eqref {Pge} holds trivially for every $n$ odd. For $n=2,4,6,10$, we have
\begin {equation}\label {10}
\begin {split}
\sum\limits_{{\bf x}\in\OL\omega_{2160}}&\sum\limits_{{\bf y}\in\OL\omega_{2160}}P_n^{(7)}({\bf x}\cdot {\bf y})=4320(P_n^{(7)}(1)+64P_n^{(7)}\(3/4\)\\
&+280P_n^{(7)}\(1/2\)+448P_n^{(7)}(1/4)+287P_n^{(7)}\(0\))=0.
\end {split}
\end {equation}
The code $\OL\omega_{2160}$ is a $7$-design (and not an $8$-design). However, it is not stiff, since $\mathcal D_4(\OL\omega_{2160})=\emptyset$ as the following lemma suggests.
\begin {lemma}\label {E_8}
For every vector ${\bf x}\in S^7$, the set $D({\bf x},\OL\omega_{2160})$ has at least five distinct elements. The only vectors ${\bf x}\in S^7$ such that $D({\bf x},\OL\omega_{2160})$ has exactly five distinct elements are those in $\OL\omega_{240}$. For each ${\bf x}\in \OL\omega_{240}$, we have $D({\bf x},\OL\omega_{2160})=\left\{0,\pm \frac {1}{2\sqrt{2}},\pm \frac {1}{\sqrt{2}}\right\}$.
\end {lemma}
\begin {proof}
Let ${\bf x}=(x_1,\ldots,x_8)\in \mathcal D_5(\OL\omega_{2160})$ be arbitrary. If non-zero coordinates of ${\bf x}$ had at least three distinct absolute values, then ${\bf x}$ would form at least $6$ distinct dot products with vectors of type II. Therefore, non-zero coordinates of ${\bf x}$ have at most two distinct absolute values. 

Assume to the contrary that non-zero coordinates of ${\bf x}$ have exactly two distinct absolute values. Denote them by $0<b<c$. Then ${\bf x}$ forms each of the dot products $\pm b,\pm c$ with vectors of type II. If ${\bf x}$ formed with some vector ${\bf z}\in \OL\omega_{2160}$ a positive dot product $u$ distinct from $b$ and $c$, since $\OL\omega_{2160}$ is antipodal, there would be a sixth dot product $-u$, contradicting the assumption that ${\bf x}\in \mathcal D_5(\OL\omega_{2160})$. Thus, $D({\bf x},\OL\omega_{2160})$ contains only two positive dot products: $b$ and $c$. 
Let $k$ coordinates of ${\bf x}$ have absolute value $b$ and $\ell$ coordinates have absolute value $c$. If it were that $\ell\geq 2$, then ${\bf x}$ would form positive dot products $\frac {2c+b+v}{2}$ and $\frac {b+v}{2}$ with appropriately chosen two vectors of type I, where $v\geq 0$ is the absolute value of one of the coordinates of ${\bf x}$. Then $\frac {2c+b+v}{2}=c$ forcing $b\leq 0$. Thus, $\ell=1$. If it were that $k\geq 3$, then ${\bf x}$ would form positive dot products $\frac {c+3b}{2}$ and $\frac {c+b}{2}$ with certain two vectors of type I. This would force $\frac {c+b}{2}=b$; that is, $c=b$. Thus, $k\leq 2$. We can now take vector ${\bf z}$ to be of type III with the coordinate with $\pm 3/4$ corresponding to a zero coordinate of ${\bf x}$ such that ${\bf x}$ forms positive dot products ${\bf x}\cdot {\bf z}=\frac {c+k b}{4}$ and $\frac {c+(k-2)b}{4}$. Then $c=\frac {c+k b}{4}\leq \frac {c+2b}{4}<\frac {3c}{4}$, which is a contradiction.

Thus, all non-zero coordinates of ${\bf x}$ have the same absolute value, which we denote by $a$. Let $n$ be the number of non-zero coordinates of ${\bf x}$. If $n=1$, then $a=1$ and ${\bf x}$ forms nine dot products, $0,\pm 1/4,\pm 1/2,\pm 3/4,\pm 1$, with points of $\OL\omega_{2160}$. If $n=3$, then ${\bf x}$ forms seven dot products, $0,\pm a/2,\pm a,\pm 3a/2$, with type I vectors. If now $4\leq n\leq 7$, then ${\bf x}$ forms nine dot products $0,\pm a/2,\pm a,\pm 3a/2,\pm 2a$ with type I vectors. Therefore, $n=2$ or $8$.

If $n=2$, then $a=1/\sqrt{2}$ and ${\bf x}\in \OL \omega_{240}$. Finally, if $n=8$, then every coordinate of ${\bf x}$ is $\pm \frac {1}{2\sqrt{2}}$. 
Assume to the contrary that ${\bf x}$ has an odd number of negative coordinates. Then for every vector ${\bf z}$ of type III, ${\bf x}\cdot {\bf z}$ is a sum of seven signed terms $\frac {1}{8\sqrt{2}}$ and one signed term $\frac {3}{8\sqrt{2}}$ with an even total number of minus signs. Then ${\bf x}\cdot {\bf z}$, in particular, has six values $\pm 2w,\pm 6w,\pm 10w$, where $w=\frac {1}{8\sqrt{2}}$. Therefore, coordinates of ${\bf x}$ have an even number of minus signs and ${\bf x}\in \OL\omega_{240}$.

Thus, if ${\bf x}\notin \OL\omega_{240}$ then ${\bf x}$ forms more than five distinct dot products with points of $\OL\omega_{2160}$. One can also verify directly that every ${\bf x}\in \OL\omega_{240}$ forms exactly five distinct dot products with points of $\OL\omega_{2160}$, which are $0,\pm \frac {1}{2\sqrt{2}},\pm \frac {1}{\sqrt{2}}$.
\end {proof}
We are now ready to state the main result of this section.

\begin {theorem}\label {2_41p}
Let $\OL\omega_{2160}=\{{\bf x}_1,\ldots,{\bf x}_{2160}\}$ be the $2_{41}$ polytope on $S^7$
and $g$ be an admissible potential function with non-negative derivatives $g^{(8)}$, $g^{(9)}$, and $g^{(10)}$ on $(-1,1)$. Then, for every point ${\bf x}^\ast\in \OL\omega_{240}$,
\begin {equation}\label {star1}
\min\limits_{{\bf x}\in S^7}\sum\limits_{i=1}^{2160}g({\bf x}\cdot{\bf x}_i)=\sum\limits_{i=1}^{2160}g({\bf x}^\ast\cdot {\bf x}_i).
\end {equation}
If, in addition, $g^{(10)}>0$ on $(-1,1)$, then the absolute minimum in \eqref {star1} is achieved only at points of the set $\OL\omega_{240}$. 
\end {theorem}
\begin {proof}
We have $\{1,2,3,4,5,6,7,9,10\}\subset\mathcal I(\OL\omega_{2160})$ in view of \eqref {10}. 
Applying Theorem \ref {generalth} with $d=7$, $m=5$, $\{t_1,\ldots,t_5\}=\left\{0,\pm \frac {1}{2\sqrt{2}},\pm \frac {1}{\sqrt{2}}\right\}$ we have $\mathcal D=\OL\omega_{240}$ and equality \eqref {star1} holds for every ${\bf x}^\ast\in\OL\omega_{240}$. If $g^{(10)}>0$ on $(-1,1)$, then, by Theorem \ref {generalth}, \eqref {star1} holds only for ${\bf x}^\ast\in \OL\omega_{240}$.
\end {proof}

\section {Three symmetrized sharp configurations on $S^{20}$ and~$S^{21}$}\label {sy}

In this section, we discuss one simple method which allows to construct some new stiff configurations and obtain their universal minima.

\begin {lemma}\label {symm}
Let $\omega_N\subset S^d$ be an $m$-stiff configuration, $m,d\geq 1$, which does not contain an antipodal pair. Let $\omega_N':=\omega_N\cup\(-\omega_N\)$ be its symmetrization. Then $\omega_N'$ is also $m$-stiff with $\mathcal D_m(\omega_N')=\mathcal D_m(\omega_N)$.
\end {lemma}
\begin {proof}
The configurations $\omega_N$ and $-\omega_N$ are disjoint and both are $(2m-1)$-designs. Then their union $\omega_N'$ has $2N$ points and is also a $(2m-1)$-design. We immediately have $\mathcal D_m(\omega_N')\subset \mathcal D_m(\omega_N)$. If ${\bf x}\in \mathcal D_m(\omega_N)$ then, by Proposition~\ref {mtight}, ${\bf x}$ forms one of the dot products $\kappa_1^m,\ldots,\kappa_m^m$ with any point from $\omega_N$. Since $\kappa_1^m,\ldots,\kappa_m^m$ are zeros of $P_m^{(d)}$, they are symmetric about the origin. Since $-{\bf x}\in \mathcal D_m(\omega_N)$, for every ${\bf y}\in -\omega_N$, we have ${\bf x}\cdot{\bf y}=-{\bf x}\cdot (-{\bf y})\in \{\kappa_1^m,\ldots,\kappa_m^m\}$ because $-{\bf y}\in \omega_N$. Then ${\bf x}\in \mathcal D_m(\omega_N')$; that is, $\mathcal D_m(\omega_N)=\mathcal D_m(\omega_N')$. Since $\mathcal D_m(\omega_N)\neq \emptyset$, we have $\mathcal D_m(\omega_N')\neq \emptyset$; that is, $\omega_N'$ is $m$-stiff.
\end {proof}

We remark that Lemmas \ref {symm} and \ref {2stiff} immediately imply Lemma \ref {dual} for $d\geq 5$ odd, since, in this case, $\OL\omega^d$ has no antipodal pair, $U_d=\OL\omega^d\cup (-\OL\omega^d)$, and $\mathcal D_2(U_d)=\omega_{2d}^\ast$. Lemma \ref {symm} also applies to the following two cases.

{\it Case I.} The Higman-Sims configuration, denoted by $\OL\omega_{100}$, is a $100$-point $3$-design on $S^{21}$, where distinct points form only dot products $-4/11$ and $1/11$ with each other (see, e.g., \cite [Table 1]{CohKum2007}). For $n=1,2,3$, we have
$$
\sum\limits_{{\bf x}\in \OL\omega_{100}}\sum\limits_{{\bf y}\in \OL\omega_{100}}P_n^{(21)}({\bf x}\cdot{\bf y})=100\(P_n^{(21)}(1)+77P_n^{(21)}\(\frac{1}{11}\)+22P_n^{(21)}\(-\frac{4}{11}\)\)=0.
$$

According to \cite {BoyDraHarSafStosharpantipodal}, this configuration is $2$-stiff. 
Paper \cite {BoyDraHarSafStosharpantipodal} finds $176$ pairs of antipodal vectors on $S^{21}$, where each vector forms only dot products $\frac {1}{\sqrt{22}}$ and $-\frac {1}{\sqrt{22}}$ with vectors from $\OL\omega_{100}$. Denote this set of $2\cdot 176=352$ vectors by $\OL\omega_{352}$.
We have $\OL\omega_{352}\subset\mathcal D_2(\OL\omega_{100})$. 
Each of these vectors is a universal minimum of $\OL\omega_{100}$ (see Theorem \ref {2m-1}). 
Since no dot product in $\OL\omega_{100}$ is $-1$, the set $\OL\omega_{100}$ does not contain an antipodal pair.

By Lemma~\ref {symm}, the symmetrized Higman-Sims configuration $\OL\omega_{200}:=\OL\omega_{100}\cup(-\OL\omega_{100})$ is $2$-stiff with $\mathcal D_2(\OL\omega_{200})=\mathcal D_2(\OL\omega_{100})\supset\OL\omega_{352}$.

{\it Case II.} Two sharp codes on $S^{20}$ can be derived from the McLaughlin configuration $\OL\omega_{275}\subset S^{21}$, which is strongly sharp (see, e.g., \cite[Table 1]{CohKum2007}). Fix a point ${\bf x}\in \OL\omega_{275}$. It forms dot product $1/6$ with $162$ points from $\OL\omega_{275}$. Let $\omega_{162}^{\bf x}$ denote the set of these $162$ points. They lie in the intersection of $S^{21}$ with the hyperplane ${\bf x}\cdot {\bf t}=1/6$. Point ${\bf x}$ forms dot product $-1/4$ with the remaining set of $112$ points from $\OL\omega_{275}$, which we denote by $\omega_{112}^{\bf x}$.

We apply homotheties to $\omega^{\bf x}_{162}$ and $\omega_{112}^{\bf x}$ to scale them to $\W S^{20}:=S^{21}\cap H$, where $H$ is the $21$-dimensional linear subspace of $\RR^{22}$ orthogonal to ${\bf x}$. Denote the resulting configurations by $\OL\omega_{162}$ and $\OL\omega_{112}$, respectively. Both $\OL\omega_{162}$ and $\OL\omega_{112}$ are $3$-designs. They are known sharp configurations, see \cite [Table 1]{CohKum2007} ($\OL\omega_{112}$ is the isomorphic subspace with $q=3$). Since any vector from $\omega_{162}^{\bf x}$ and any vector from $\omega^{\bf x}_{112}$ form only dot products $1/6$ or $-1/4$ with each other, any vector from $\OL\omega_{162}$ and any vector from $\OL\omega_{112}$ form only dot products $\frac {1}{\sqrt{21}}$ or $-\frac {1}{\sqrt{21}}$. Then both configurations are $2$-stiff with $\OL\omega_{112}\subset\mathcal D_2(\OL\omega_{162})$ and $\OL\omega_{162}\subset\mathcal D_2(\OL\omega_{112})$ (this was observed in \cite {BoyDraHarSafStosharpantipodal}). 

 Since any two distinct points in $\omega_{162}^{\bf x}$ form dot products $1/6$ or $-1/4$ with each other, any two distinct points in $\OL\omega_{162}$ form dot products $1/7$ or $-2/7$. Since any two distinct points in $\omega_{112}^{\bf x}$ form dot products $1/6$ or $-1/4$ with each other, any two distinct points in $\OL\omega_{112}$ form dot products $1/9$ or $-1/3$. 
In particular, $\OL\omega_{162}$ and $\OL\omega_{112}$ do not contain an antipodal pair. Let $\OL\omega_{324}:=\OL\omega_{162}\cup (-\OL\omega_{162})$ and $\OL\omega_{224}:=\OL\omega_{112}\cup (-\OL\omega_{112})$ be their symmetrization about the origin. It is not difficult to see that $\OL\omega_{224}\subset\mathcal D_2(\OL\omega_{162})$ and $\OL\omega_{324}\subset\mathcal D_2(\OL\omega_{112})$.
Each point of $\OL\omega_{224}$ is a universal minimum of $\OL\omega_{162}$ and each point of $\OL\omega_{324}$ is a universal minimum of $\OL\omega_{112}$ (see Theorem \ref {2m-1}).

By Lemma~\ref {symm}, both $\OL\omega_{324}$ and $\OL\omega_{224}$ are $2$-stiff with
$\mathcal D_2(\OL\omega_{324})=\mathcal D_2(\OL\omega_{162})\supset\OL\omega_{224}$ and $\mathcal D_2(\OL\omega_{224})=\mathcal D_2(\OL\omega_{112})\supset\OL\omega_{324}$.

Concluding paragraphs in Cases I and II and Theorem \ref {2m-1} imply the following.
\begin {proposition}\label {thsym}
Let $g$ be an admissible potential function with a convex derivative $g''$ on $(-1,1)$. Then
\begin {itemize}
\item[(i)]
every point of the configuration $\OL\omega_{352}$ is a point of absolute minimum over $S^{21}$ of the potential $p^g(\cdot,\OL\omega_{200})$ of the symmetrized Higman-Sims configuration~$\OL\omega_{200}$;
\item [(ii)]
every point of the configuration $\OL\omega_{224}$ is a point of absolute minimum over $\W S^{20}$ of the potential $p^g(\cdot,\OL\omega_{324})$;
\item[(iii)]
every point of the configuration $\OL\omega_{324}$ is a point of absolute minimum over $\W S^{20}$ of the potential $p^g(\cdot,\OL\omega_{224})$.
\end {itemize}
\end {proposition}

\section {Certain properties of general stiff configurations}\label {scstiff}

Every time we have a stiff configuration, in view of Theorem \ref {2m-1}, we automatically have its universal minima (the dual configuration). Moreover, every stiff configuration attains the Fazekas-Levenshtein bound for covering \cite [Theorem~2]{FazLev1995}. Therefore, it is important to study stiff codes and their duals in general.
In this section, we characterize $1$-stiff configurations on $S^d$, $m$-stiff configurations on $S^1$, and their duals and also prove some basic properties of stiff configurations and of their duals.
%
We call the point ${\bf c}=\frac {1}{N}\sum_{i=1}^{N}{\bf x}_i$ {\it the center of mass} of a configuration $\omega_N=\{{\bf x}_1,\ldots,{\bf x}_N\}$.
\begin {proposition}\label {1stiff}
Let $d\geq 1$. A configuration $\omega_N\subset S^d$, $N\geq 1$, is $1$-stiff if and only if its center of mass is at the origin and $\omega_N$ is contained in a $d$-dimensional linear subspace of $\RR^{d+1}$.
\end {proposition}
\begin {proof}
The proposition follows from the fact that a point configuration is a spherical $1$-design if and only if its center of mass is located at the origin and the fact that a hyperplane containing $\omega_N$ also contains its center of mass.
\end {proof}

We next describe the dual of a $1$-stiff configuration.
\begin {proposition}\label {perp}
Let $\omega_N\subset S^d$, $d\geq 1$, be a $1$-stiff configuration. Then $\mathcal D_1(\omega_N)=L^\bot\cap S^d$, where $L$ is the linear subspace of $\RR^{d+1}$ spanned by $\omega_N$. If $k:={\rm dim}\ \! L\leq d-1$, then $\mathcal D_1(\omega_N)$ is a sphere in a $(d+1-k)$-dimesnional subspace of $\RR^{d+1}$. If $k=d$, then $\mathcal D_1(\omega_N)=\{{\bf a},-{\bf a}\}$ for some ${\bf a}\in S^d$, which is a $1$-stiff configuration.
\end {proposition}
\begin {proof}
Let ${\bf z}$ be any vector in $\mathcal D_1(\omega_N)$. By Proposition \ref {mtight}, we have ${\bf z}\cdot {\bf y}=0$ for every ${\bf y}\in \omega_N$, since $0$ is the only root of $P_1^{(d)}$. Then ${\bf z}\bot L$; i.e., ${\bf z}\in L^\bot\cap S^d$. If ${\bf z}$ is any vector in $L^\bot\cap S^d$, then it forms only one dot product (which is $0$) with any point from $\omega_N$; that is, ${\bf z}\in \mathcal D_1(\omega_N)$. The rest of Proposition \ref {perp} follows immediately.
%
%
%
\end {proof}

We now charactirize stiff configurations on $S^1$.
\begin {proposition}\label {S1}
For every $m\geq 1$, a configuration on $S^1$ is $m$-stiff  if and only if it is a regular $2m$-gon.
\end {proposition}
We remark that the regular $2m$-gon $\W\omega_{2m}$ on $S^1$ is antipodal and its dual $\mathcal D_m(\W\omega_{2m})$ is another regular $2m$-gon with $\mathcal D_m(\mathcal D_m(\W\omega_{2m}))=\W\omega_{2m}$. 
\begin {proof}[Proof of Proposition \ref {S1}]
If $\omega_N=\W \omega_{2m}$ then it is a $(2m-1)$-design and the midpoint ${\bf y}$ of the arc joining any two neighboring vertices forms $m$ distinct values of dot product with points from $\omega_N$; i.e., $\omega_N$ is $m$-stiff.

By Proposition \ref {mtight}, the point ${\bf y}$ forms dot products $\kappa_1^m,\ldots,\kappa_m^m$ (zeros of $P_m^{(1)}$) with points of $\W\omega_{2m}$, each with frequency $2=2ma_0(\varphi_i)$, $i=1,\ldots,m$, where $\{\varphi_1,\ldots,\varphi_m\}$ is the fundamental system of polynomials for the nodes $\kappa_i^m$. We have
\begin {equation}\label {phi1}
a_0(\varphi_i)=\frac{1}{m} \ \ \text{and}\ \ \kappa_i^m=\cos \frac {(2i-1)\pi}{2m}, \ \ \ \ i=1,\ldots,m.
\end {equation}
Assume that $\omega_N$ is $m$-stiff.
Let ${\bf z}$ be any point in $\mathcal D_m(\omega_N)$. By Proposition~\ref {mtight}, point ${\bf z}$ forms dot products $\kappa_1^m,\ldots,\kappa_m^m$ with points of $\omega_N$, where $\kappa^m_i$, $i=1,\ldots,m$, are the zeros of $P_m^{(1)}$. Then $\omega_N$ is contained in the set of points of intersection of $S^1$ with $m$ parallel lines; that is, $\#\omega_N\leq 2m$. By Proposition \ref {mtight} and \eqref {phi1}, the frequency $M_i$ of the dot product $\kappa^m_i$ is $M_i=Na_0(\varphi_i)=N/m\leq 2$, $i=1,\ldots,m$. Hence, frequencies $M_i$ are equal and each of the $m$ parallel lines contains the same number of points from $\omega_N$ (one or two). Assume to the contrary that each $M_i$ equals $1$. Then $\omega_N$ has only $m$ points, which means that any point ${\bf y}\in\omega_N$ forms at most $m$ distinct values of the dot product with points of $\omega_N$; i.e., ${\bf y}\in \mathcal D_m(\omega_N)$. One of dot products is~$1$, while by Proposition \ref {mtight}, these dot products must be zeros $\kappa_i^m$ of $P_m^{(1)}$ none of which is $1$. This contradiction shows that $M_i=2$, $i=1,\ldots,m$, and, hence $\#\omega_N=2m$. Vector ${\bf z}$ forms each of the angles $\frac {(2i-1)\pi}{2m}=\arccos \kappa_i^m$ with exactly two points from $\omega_N$, $i=1,\ldots,m$. Then $\omega_N$ is a regular $2m$-gon.
\end {proof}

We next prove the following basic statement.
\begin {proposition}\label {symm'}
For any configuration $\omega_N\subset S^d$, $d\geq 1$, with $\mathcal D_m(\omega_N)\neq \emptyset$, $m\geq 1$,
the set $\mathcal D_m(\omega_N)$ is antipodal. 
If $\omega_N$ is $m$-stiff, then 
$\omega_N\subset \mathcal D_m(\mathcal D_m(\omega_N))$. 
If, in addition, $\mathcal D_m(\omega_N)$ is finite and $m$-stiff, then $\mathcal D_m(\mathcal D_m(\mathcal D_m(\omega_N)))=\mathcal D_m(\omega_N)$.
\end {proposition}
Though the dual $\mathcal D_m(\omega_N)$ is antipodal, this is not always true for the $m$-stiff configuration $\omega_N$ itself. For example, the $d$-demicube $\OL\omega^d$ on $S^{d-1}$ is not antipodal for any $d\geq 3$ odd. 

In view of Proposition \ref {symm'}, equality $\omega_N=\mathcal D_m(\mathcal D_m(\omega_N))$ with an $m$-stiff $\omega_N$ implies that $\omega_N$ is antipodal. However, the inclusion $\omega_N\subset \mathcal D_m(\mathcal D_m(\omega_N))$ can be sometimes strict (even when both $\omega_N$ and $\mathcal D_m(\omega_N)$ are antipodal and $m$-stiff). This is the case, for example, for $\omega_N=\OL\omega^d$ and any $d\geq 5$ in view of Lemma~\ref {dual} and the fact that the dual of the cross-polytope $\omega_{2d}^\ast$ is the whole cube $U_d$ (for $d\geq 5$ even, both $\OL\omega^d$ and $\mathcal D_2(\OL\omega^d)=\omega^\ast_{2d}$ are antipodal and $2$-stiff).
At the same time, every stiff configuration from \cite [Table 3]{Borstiff} coincides with the dual of its dual.

Some other examples of non-antipodal stiff configurations $\omega_N$ are given in \cite {BoyDraHarSafStosharpantipodal}. Since, by Proposition \ref {symm'}, the duals of their duals are antipodal, the inclusion $\omega_N\subset \mathcal D_m(\mathcal D_m(\omega_N))$ is strict as well. 
\begin {proof}[Proof of Proposition \ref {symm'}]
Let $\omega_N\subset S^d$ be arbitrary with $\mathcal D_m(\omega_N)\neq \emptyset$.
For any point ${\bf z}\in \mathcal D_m(\omega_N)$, the point $-{\bf z}$ also forms at most $m$ distinct dot products with points of $\omega_N$; that is $-{\bf z}\in\mathcal D_m(\omega_N)$ and the set $\mathcal D_m(\omega_N)$ is antipodal. Choose any point ${\bf z}$ in an $m$-stiff $\omega_N$. For any point ${\bf y}\in \mathcal D_m(\omega_N)$, by Proposition \ref {mtight}, we have ${\bf y}\cdot {\bf z}\in \{\kappa_1^m,\ldots,\kappa_m^m\}$; that is, ${\bf z}$ forms at most $m$ distinct dot products with points of $\mathcal D_m(\omega_N)$. Then ${\bf z}\in \mathcal D_m(\mathcal D_m(\omega_N))$ and $\omega_N\subset\mathcal D_m(\mathcal D_m(\omega_N))$.
Assume additionally that $X:=\mathcal D_m(\omega_N)$ is finite and $m$-stiff. The inclusion $\omega_N\subset \mathcal D_m(\mathcal D_m(\omega_N))$ implies that $X=\mathcal D_m(\omega_N)\supset \mathcal D_m(\mathcal D_m(X))$. Since $X$ is $m$-stiff, we have the opposite inclusion.
\end {proof}

 We say that a point set in $\RR^{d+1}$ is {\it in general position} if it is not contained in any hyperplane. One can construct plenty of examples of $m$-stiff configurations, $m\geq 2$, whose dual is not in general position. For instance, start with the cube $U_3$ inscribed in $S^2$ and let $\alpha_1$ and $\alpha_2$ be the parallel planes containing some two paralel facets of $U_3$. For a given $n\geq 2$, we rotate the cube $U_3$ about the axis $\ell 
$ perpendicular to the planes $\alpha_1$ and $\alpha_2$ and passing through the origin by angles $\frac {\pi k}{2n}$, $k=0,1,\ldots,n-1$, and let $\omega_N$ be the union of the resulting $n$ cubes. Then $\omega_N$ is a $3$-design as a disjoint union of finitely many $3$-designs. Since $\omega_N$ is still contained in planes $\alpha_1$ and $\alpha_2$, it is $2$-stiff. However, its dual is $\mathcal D_2(\omega_N)=\{{\bf a},-{\bf a}\}$, where ${\bf a}$ is a unit vector parallel to the axis $\ell$. The dual is not in general position. It is also only $1$-stiff. This example can, of course, be extended to other dimensions and other initial configurations.

\begin {proposition}\label {gp}
Let $\omega_N\subset S^d$, $d\geq 1$, be an $m$-stiff configuration, $m\geq 2$. Then $\mathcal D_m(\omega_N)$ contains at most $m^{d+1}$ points. If $\mathcal D_m(\omega_N)$ is not in general position, then $\mathcal D_m(\omega_N)$ is $1$-stiff and, hence, not $m$-stiff.
\end {proposition}
Proposition \ref {gp} implies that an $m$-stiff configuration on $S^d$, $m\geq 2$, cannot have more than $m^{d+1}$ universal minima.
This cardinality bound can be achieved: take $\omega_N$ to be the set of vertices of a regular cross-polytope inscribed in $S^d$. It is $2$-stiff and its dual is a cube inscribed in $S^d$, which has exactly $2^{d+1}$ vertices.
At the same time, no upper bound depending only on $m$ and $d$ can be written for the cardinality of an $m$-stiff configuration $\omega_N$ itself, $d\geq 2$ (when $\omega_N$ exists for those $m$ and $d$), see Proposition \ref {semigroup} below.

Another interesting question related to Proposition \ref {gp} is about the general assumptions under which the dual of a given $m$-stiff configuration, $m\geq 2$, is in general position and whether this is sufficient for the dual to be also $m$-stiff. The dual is $m$-stiff with the same $m$, for example, for every stiff configuration mentioned in \cite [Table 3]{Borstiff} and for the $d$-demicube, $d\geq 4$. 

\begin {proof}[Proof of Proposition \ref {gp}]
Since $\omega_N$ is at least a $3$-design, it is in general position. If it were not, then $\omega_N$ would not be a $2$-design: for the polynomial $p({\bf x})=({\bf x}\cdot {\bf v}-\alpha)^2$, where ${\bf x}\cdot {\bf v}=\alpha$ is an equation of the hyperplane containing $\omega_N$, the average of $p$ over $\omega_N$ is zero while the average of $p$ over $S^d$ is positive. Since $\omega_N$ is in general position, it contains a linearly independent subset $\{{\bf y}_1,\ldots,{\bf y}_{d+1}\}$. 

Let ${\bf z}$ be any point in $\mathcal D_m(\omega_N)$. Then, by Proposition \ref {mtight}, we have ${\bf z}\cdot {\bf y}_i\in K_m:=\{\kappa_1^m,\ldots,\kappa_m^m\}$, $i=1,\ldots,d+1$, where $\kappa_j^m$'s are zeros of $P_m^{(d)}$. Consequently, ${\bf z}$ is a solution to a linear system ${\bf y}_j\cdot {\bf z}=\alpha_j$, $j=1,\ldots,d+1$, where $\alpha_1,\ldots,\alpha_{d+1}\in K_m$. Since this system has a non-singular coefficient matrix, it has a unique solution, which may or may not be on $S^d$. Since there are $m^{d+1}$ possible vectors of right-hand sides for this system, there are at most $m^{d+1}$ points in $\mathcal D_m(\omega_N)$.

Assume that $\mathcal D_m(\omega_N)$ is contained in some hyperplane $H$. By Proposition~\ref {symm'}, it is antipodal. Then its center of mass is at the origin; that is $H$ is a $d$-dimensional linear subspace of $\RR^{d+1}$. By Proposition \ref {1stiff}, $\mathcal D_m(\omega_N)$ is $1$-stiff. Since $\mathcal D_m(\omega_N)$ is contained in one hyperplane, by the above argument, it cannot be a $2$-design.
Then $\mathcal D_m(\omega_N)$ cannot be $m$-stiff.
\end {proof}

The statement below implies that the set of possible cardinalities of $m$-stiff configurations on $S^d$, $d\geq 2$ (provided that an $m$-stiff configuration exists on $S^d$), forms an additive semigroup (in particular, it is not bounded above). This is not the case for $d=1$ in view of Proposition~\ref {S1}. In the case $m=1$ and $d\geq 2$, this semigroup is the set of all integers $N\geq 2$, see Proposition \ref {1stiff}.
\begin {proposition}\label {semigroup}
Suppose that for given $m\geq 1$ and $d\geq 2$, there exist $m$-stiff configurations on $S^d$ of cardinalities $N_1$ and $N_2$. Then there is an $m$-stiff configuration on $S^d$ of cardinality $N_1+N_2$.
\end {proposition}
Proposition \ref {semigroup} and B\'ezout's identity imply that if, for a given pair $(m,d)$, $d\geq 2$, cardinalities of some two $m$-stiff configurations on $S^d$ have greatest common divisor $\delta$, then for any sufficiently large multiple $N$ of $\delta$, there exists an $m$-stiff configuration on $S^d$ of cardinality $N$.
\begin {proof}
Let $\omega_{N_i}\subset S^d$ be an $m$-stiff configuration of cardinality $N_i$, $i=1,2$, and let ${\bf z}_i\in \mathcal D_m(\omega_{N_i})$, $i=1,2$, be chosen so that ${\bf z}_1\neq {\bf z}_2$ (if it happens that ${\bf z}_1={\bf z}_2$, we choose $-{\bf z}_2$ instead of ${\bf z}_2$). We will construct an $m$-stiff configuration of cardinality $N_1+N_2$. By Proposition \ref {mtight}, both vectors ${\bf z}_1$ and ${\bf z}_2$ form only dot products $\kappa_1^m,\ldots,\kappa_m^m$ with points from the corresponding configuration $\omega_{N_i}$. Let $H$ be the $d$-dimensional subspace of $\RR^{d+1}$, which is the perpendicular bisector for the line segment $[{\bf z}_1,{\bf z}_2]$. Let $r_H:\RR^{d+1}\to\RR^{d+1}$ be the reflection transformation about the subspace $H$ and let $U$ be its matrix ($U$ is orthogonal). Then for every ${\bf x}\in \omega_{N_1}':=r_H(\omega_{N_1})$, there is ${\bf y}\in \omega_{N_1}$ such that ${\bf x}=U{\bf y}$ and 
\begin {equation}\label {z2}
{\bf x}\cdot {\bf z}_2=U{\bf y}\cdot {\bf z}_2={\bf y}\cdot U^T{\bf z}_2={\bf y}\cdot {\bf z}_1=\kappa_i^m \ \ \ \text{for some}\ \ \ i=1,\ldots,m.
\end {equation}
The set $\omega_{N_1}'$ is $m$-stiff.
Denote by ${\bf a}\in S^d$ a vector perpendicular to ${\bf z}_2$, not contained in any subspace ${\rm span}\{{\bf w}-{\bf v}\}$, where ${\bf w}\in \omega_{N_2}$ and ${\bf v}\in \omega_{N_1}'$, and ${\bf a}$ is not perpendicular to any vector from $\omega_{N_2}$. Such a vector ${\bf a}$ exists, since all these conditions delete a set of $(d-1)$-dimensional measure zero from $S^d\cap \{{\bf z}_2\}^\bot$. 

Then $L=\{{\bf a}\}^\bot$ is disjoint with $\omega_{N_2}$. Furthermore, ${\bf z}_2\in L$ and the configuration $\omega_{N_1}'':=r_L(\omega_{N_1}')$ is disjoint with $\omega_{N_2}$. If it were not, then $L$ would be the perpendicular bisector for a line segment whose one endpoint ${\bf w}_1$ is in $\omega_{N_2}$ and the other endpoint ${\bf v}_1$ is in $\omega_{N_1}'$. We have ${\bf w}_1\neq{\bf v}_1$, since, otherwise, ${\bf w}_1\in L$ and, hence, ${\bf a}\bot {\bf w}_1$. Then ${\bf w}_1-{\bf v}_1$ is a non-zero vector perpendicular to $L$ and ${\bf a}\in {\rm span}\{{\bf w}_1-{\bf v}_1\}$ contradicting the choice of ${\bf a}$. 

Let $V$ be the matrix of the reflection transformation $r_L$. Then for every ${\bf z}\in \omega_{N_1}''$, there is ${\bf x}\in \omega_{N_1}'$ such that ${\bf z}=V{\bf x}$ and using \eqref {z2} we have
$$
{\bf z}\cdot {\bf z}_2=V{\bf x}\cdot {\bf z}_2={\bf x}\cdot V^T{\bf z}_2={\bf x}\cdot {\bf z}_2=\kappa_i^m\ \ \ \text{for some}\ \ \ i=1,\ldots,m.
$$
Thus, ${\bf z}_2\in \mathcal D_m(\omega_{N_1}''\cup\omega_{N_2})$ and $\omega_{N_1}''\cup\omega_{N_2}$ is a disjoint union of two $(2m-1)$-designs. Then it is also a $(2m-1)$-design, and, hence is an $m$-stiff configuration of cardinality $N_1+N_2$.
\end {proof}

}


\begin {thebibliography}{99}
\bibitem {Bil2015}
M. Bilogliadov, Equilibria of Riesz potentials generated by point charges at the roots of unity, {\it Comput. Methods Funct. Theory} {\bf 15} (2015), no. 4, 471--491.
\bibitem {Borstiff}
S.V. Borodachov, Absolute minima of potentials of a certain class of spherical designs (submitted), https://arxiv.org/abs/2212.04594.
\bibitem {Borsymmetric}
S.V. Borodachov, Absolute minima of potentials of certain regular spherical configurations (submitted), https://arxiv.org/abs/2210.04295.
\bibitem {BorBatumi}
S.V. Borodachov, Extreme values of potentials of spherical designs and the polarization problem, {\it XII Annual International Conference of the Georgian Mathematical Union}, Batumi State University, Georgia, August 29--September 3.
\bibitem {BorMinMax}
S.V. Borodachov, Min-max polarization for certain classes of sharp configurations on the sphere (submitted), https://arxiv.org/abs/2203.13756.
\bibitem{Bor2022talk}
S.V. Borodachov, Min-max polarization for certain classes of sharp configurations on the sphere, {\it Workshop "Optimal Point Configurations on Manifolds"}, ESI, Vienna, January 17--21, 2022. https://www.youtube.com/watch?v=L-szPTFMsX8
\bibitem {Borsimplex}
S.V. Borodachov, Polarization problem on a higher-dimensional sphere for a simplex, {\it Discrete and Computational Geometry} {\bf 67} (2022), 525--542.

\bibitem{BorHarSafbook}
S. Borodachov, D. Hardin, E. Saff, {\it Discrete Energy on Rectifiable Sets}. Springer, 2019.
\bibitem {BoyDraHarSafSto600cell}
P.G. Boyvalenkov, P.D. Dragnev, D.P. Hardin, E.B. Saff, M.M. Stoyanova, On polarization of spherical codes and designs (submitted), https://arxiv.org/abs/2207.08807.
\bibitem {BoyDraHarSafStosharpantipodal}
P.G. Boyvalenkov, P.D. Dragnev, D.P. Hardin, E.B. Saff, M.M. Stoyanova, Universal minima of discrete potentials for sharp spherical codes, https://arxiv.org/pdf/2211.00092.pdf.
\bibitem{CohKum2007}
H. Cohn, A. Kumar, Universally optimal distribution of points on spheres, {\it J. Amer. Math. Soc.} {\bf 20} (2007), no. 1, 99--148.
\bibitem {DelGoeSei1977}
P. Delsarte, J.M. Goethals, J.J. Seidel, Spherical codes and designs, {\it Geometriae Dedicata}, {\bf 6} (1977), no. 3, 363--388.
\bibitem {FazLev1995}
G. Fazekas, V.I. Levenshtein, On upper bounds for code distance and covering radius of designs in polynomial metric spaces, {\it Journal of Combinatorial Theory, Series A} {\bf 70} (1995), no. 2, 267--288.
\bibitem {GioKhi2020}
G. Giorgadze, G. Khimshiashvili, Stable equilibria of three constrained unit charges, {\it Proc. I. Vekua Inst. Appl. Math.} {\bf 70} (2020), 25--31.
\bibitem {Gos1900}
T. Gosset, On the regular and semi-regular figures in space of $n$ dimensions, {\it Messenger of Mathematics}, Macmillan, 1900.
\bibitem {HarKenSaf2013}
D. Hardin, A. Kendall, E. Saff,
Polarization optimality of equally spaced points on the circle for discrete potentials.
{\it Discrete Comput. Geom.} {\bf 50} (2013), no. 1, 236--243. 
\bibitem {Lev1979}
V.I.Levenshtein, On bounds for packings in $n$-dimensional Euclidean space,
{\it Soviet Math. Dokladi} {\bf 20} (1979), 417--421.
\bibitem {Lev1992}
V.I. Levenshtein, Designs as maximum codes in polynomial metric spaces,
{\it Acta Appl. Math.} {\bf 25} (1992), 1--82.
\bibitem {Lev1998}
V.I. Levenshtein, Universal bounds for codes and designs, Chapter 6 in
{\it Handbook of Coding Theory}, V. Pless and W.C. Huffman, Eds., Elsevier Science B.V., 1998.
\bibitem {NikRaf2011}
N. Nikolov, R. Rafailov, On the sum of powered distances to certain sets of points on the circle, {\it Pacific J. Math.} {\bf 253} (2011), no. 1, 157--168.
\bibitem{NikRaf2013}
N. Nikolov, R. Rafailov, On extremums of sums of powered distances to a finite set of points. {\it Geom. Dedicata} {\bf 167} (2013), 69--89.
\bibitem {Sto1975circle}
K. Stolarsky, The sum of the distances to certain pointsets on the unit circle, {\it Pacific J. Math.} {\bf 59} (1975), no. 1, 241--251. 
\bibitem {Sto1975}
K. Stolarsky, The sum of the distances to $N$ points on a sphere, {\it Pacific J. Math.} {\bf 57} (1975), no. 2, 563--573.
\bibitem {Yud1992}
V.A. Yudin, Minimum potential energy of a point system of charges, {\it Diskrete Math. Appl.} {\bf 3} (1993), no. 1, 75--81.
\bibitem {Sze1975}
G. Szeg\"o, {\it Orthogonal polynomials}. Fourth edition. American Mathematical Society, Colloquium Publications, Vol. XXIII. American Mathematical Society, Providence, R.I., 1975.
\end {thebibliography}
\end {document}